\documentclass{amsproc}
\usepackage{amsmath}
\usepackage{amsfonts}
\usepackage{amssymb}
\usepackage{epsfig,latexsym,graphicx,stmaryrd}

\makeindex

\newtheorem{theorem}{Theorem}[section]
\newtheorem{lemma}[theorem]{Lemma}
\newtheorem{corollary}[theorem]{Corollary}

\theoremstyle{definition}
\newtheorem{definition}[theorem]{Definition}

\theoremstyle{remark}
\newtheorem{remark}[theorem]{Remark}

\numberwithin{equation}{section}

\newcommand{\ignore}[1]{}

\newcommand{\ip}[2]{{\langle#1,#2\rangle}}
\newcommand{\ipp}[2]{{\langle\langle#1,#2\rangle\rangle}}

\newcommand{\norm}[1]{{\|#1\|}}
\newcommand{\outp}[2]{{\llbracket #1,#2 \rrbracket }}

\def \cal{\mathbb}

\def \beq{\begin{eqnarray*}}
\def \eeq{\end{eqnarray*}}

\newcommand{\N}{\mathbb N}
\def \R{\mathbb{R}}
\def \C{\mathbb{C}}
\def \Z{\mathbb{Z}}

\def \PP{\mathbb{P}}
\def \E{\mathbb{\cal E}}

\newcommand{\fc}{{\mathcal{F}}}
\newcommand{\I}{{\mathbb I}}
\newcommand{\Ac}{\mathbb  A}
\newcommand{\Vc}{{\mathcal{V}}}

\newcommand{\Acr}{{\mathcal{A}}}
\newcommand{\Nc}{{\mathcal{N}}}
\newcommand{\Pp}{{\mathcal P}}

\newcommand{\Rbc}{R}
\newcommand{\Rsc}{\mathcal{R}}
\newcommand{\Zc}{\mathcal{Z}}
\newcommand{\Ssc}{\mathcal{S}}

\newcommand{\PPRR}{{\Pi_{z_0}}}

\renewcommand{\SS}{\mathcal{S}}
\newcommand{\Sonezero}{{\SS^{1,0}}}
\newcommand{\Soneone}{{\SS^{1,1}}}

\newcommand{\Ho}{\mathring{\widehat{\C^n}}}
\newcommand{\Vo}{\mathring{\hat{V}}}

\newcommand{\eps}{{\varepsilon}}

\newcommand{\Sym}{\operatorname{Sym}}
\newcommand{\trace}{\operatorname{trace}}
\newcommand{\real}{\operatorname{real}}
\newcommand{\imag}{\operatorname{imag}}
\newcommand{\Lip}{\operatorname{Lip}}
\newcommand{\supp}{\operatorname{supp}}
\renewcommand{\span}{\operatorname{span}}  
\newcommand{\Cov}{\operatorname{Cov}}
\newcommand{\rank}{\operatorname{rank}}
\newcommand{\Card}{\operatorname{Card}}
\newcommand{\snr}{\operatorname{snr}}

\newcommand{\AWGN}{\textnormal{AWGN}}
\newcommand{\nonAWGN}{\textnormal{nonAWGN}}
\newcommand{\est}{\textnormal{est}}

\begin{document}

\title[Phase Retrieval]{Frames and Phaseless Reconstruction}


\author{Radu Balan}
\address{Department of Mathematics, Center for Scientific Computation and Mathematical Modeling, Norbert Wiener Center \\ University of Maryland, College Park MD 20742}
\email{rvbalan@math.umd.edu} 

\thanks{Financial support from NSF Grants DMS-1109498 and DMS-1413249 is gratefully acknowledged.}


\date{}

\begin{abstract}
Frame design for phaseless reconstruction is now part of the broader
problem of nonlinear reconstruction and is an emerging topic in
harmonic analysis. The problem of phaseless reconstruction can be
simply stated as follows. Given the magnitudes of the coefficients generated by
 a linear redundant system (frame), we want to reconstruct
the unknown input. This problem first occurred in X-ray
crystallography starting in the early 20th century. 
The same nonlinear reconstruction problem shows up in
speech processing, particularly in speech recognition.

In this lecture we shall cover existing analysis results as well as
stability bounds for signal recovery including:
necessary and sufficient conditions for injectivity,
Lipschitz bounds of the nonlinear map and its left inverses,
stochastic performance bounds, and algorithms for signal recovery.
\end{abstract}

\maketitle

\section{Introduction}
These lecture notes concern the problem of finite dimensional vector reconstruction from magnitudes of frame coefficients.

Variants of this problem appear in several areas of engineering and science. In particular in X-ray crystallography one measures the magnitudes of the Fourier transform of the electron density from which one infers the atomic structure of the crystal \cite{fienup82}. 
In speech processing, automatic speech recognition engines typically use cepstral coefficients, which are absolute values of linear combinations of the short time Fourier transform coefficients \cite{HLO80,Bal2010}. The \index{reconstruction! phaseless}phaseless reconstruction problem unifies these and other similar problems \cite{BCE06}.

While the problem can be stated in the more general context of infinite dimensional Hilbert spaces, in these lectures we focus exclusively on the finite dimensional case.
In this case any spanning set is a frame (see \cite{Cass-artofframe} for a complete definition and list of properties).
Specifically let $H=\C^n$ denote the $n$ dimensional complex Hilbert space and let 
$\fc=\{f_1,\ldots,f_m\}$ be a set of $m\geq n$ vectors that span $H$. 
Fix a real linear space $V$,\index{space! real linear} that is also a subset of $H$, $V\subset H$.
Our problem is to study when a vector $x\in V$ can be reconstructed from
magnitudes of its frame coefficients  $\{|\ip{x}{f_k}|\, ,\,1\leq k\leq m\}$  and how to do so efficiently. 
This setup covers both the real case and the complex case as studied before in literature: in the real case $\fc\subset V=\R^n$; in the complex case $V=H=\C^n$.
Note we assume $V$ is a real linear space which may not be closed under multiplication with complex scalars.
While the analysis both of the real case and complex case is presented in a unified way, the reader should be aware that the two cases are not equally easy.
In particular the geometric criterion (4) of Theorem \ref{th3.2} does not have a counterpart in the complex case.  

Consider the following additional notation. Let 
\begin{equation}
T:H\rightarrow \C^m~~,~~(T(x))_{k}=\ip{x}{f_k}~,~1\leq k\leq m
\end{equation}
denote the frame {\em analysis map}. Its adjoint is called the {\em synthesis map} and is defined by
\begin{equation}
T^*:\C^m\rightarrow H~~,~~T^*(c)=\sum_{k=1}^m c_k f_k .
\end{equation}
We now define the main nonlinear function we discuss in this paper 
$$x\mapsto (|\ip{x}{f_k}|)_{1\leq k\leq m}.$$
For two vectors $x,y\in H$, consider the equivalence relation $x\sim y$ if and only if there is a constant $c$ of magnitude $1$ so that $x=cy$.
Thus $x\sim y$ if and only if $x=e^{i\varphi}y$ for some real $\varphi$. 
Let $\hat{H}=H/\sim$ denote the quotient space. Note the nonlinear map is well defined on $\hat{H}$ since 
$|\ip{cx}{f_k}|=|\ip{x}{f_k}|$ for all scalars $c$
with $|c|=1$. We let $\alpha$ denote the quotient map
\begin{equation}
\alpha:\hat{H}\rightarrow \R^m~~,~~(\alpha(x))_k = |\ip{x}{f_k}|~,~1\leq k\leq m .
\end{equation}
For purposes that will become clear later, let us also define the map
\begin{equation}
\beta:\hat{H}\rightarrow \R^m~~,~~(\beta(x))_k = |\ip{x}{f_k}|^2~,~1\leq k\leq m .
\end{equation}

For the subspace $V$ denote by $\hat{V}$ the set of equivalence classes $\hat{V}=\{ \hat{x}~,~x\in V\}$.

\begin{definition}
The frame $\fc$ is called a {\em phase retrievable frame\index{frame! phase retrievable frame} with respect to a set $V$} if the restriction
 $\alpha{|}_{\hat{V}}$ is injective.
\end{definition}

In these lecture notes we study the following problems:
\begin{enumerate}
\item Find necessary and sufficient conditions for $\alpha{|}_{\hat{V}}$ to be a one-to-one (injective) map;
\item Study Lipschitz properties of the maps $\alpha$, $\beta$ and their inverses;
\item Study robustness guarantees (such as Cramer-Rao Lower Bounds) for any inversion algorithm;
\item Recovery using convex algorithms (Linear Tensor recovery, and PhaseLift);
\item Recovery using iterative algorithms (Gerchberg-Saxton, Wirtinger flow, regularized least-squares).
\end{enumerate}

\section{Geometry of $\hat{H}$ and $\SS^{p,q}$ Spaces}

\subsection{$\hat{H}$} 
Recall $\hat{H}=\widehat{\C^n}=\C^n/\sim=\C^n/T^1$ where $T^1=\{z\in \C~,~|z|=1\}$.
Algebraically $\widehat{\C^n}$ is a homogeneous space invariant to multiplications by positive real scalars. In particular any $x\in\widehat{\C^n}\setminus\{0\}$ has a unique decomposition
$x=rp$, where $r=\norm{x}>0$ and $p\in\C\PP^{n-1}$ is in the projective space\index{space! projective space} $\C\PP^{n-1}=\PP(\C^n)$. Thus, topologically,
\[ \widehat{\C^n} = \{0\} \cup \left( (0,\infty)\times \C\PP^{n-1} \right). \]
The subset
\[ \Ho = \widehat{\C^n}\setminus \{0\}  = (0,\infty) \times \C\PP^{n-1} \]
is a real analytic manifold. 

Now consider the set $\hat{V}$ of equivalence classes associated to vectors in $V$. Similar to $\hat{H}$,
$\hat{V}$ admits the following decomposition
\[ \hat{V} = \{0\} \cup \left( (0,\infty) \times \PP(V) \right), \]
where $\PP(V)=\{~\{zx~,~z\in \C\}~,~x\in V,x\neq 0\}$ denotes the {\em projective space associated to $V$}.
The interior subset
\[ \Vo = \hat{V}\setminus \{0\}  = (0,\infty) \times \PP(V) \]
is a real analytic manifold of (real) dimension $1+\dim_\R\PP(V)$.

Two important cases are as follows:
\begin{itemize}
\item Real case. $V=\R^n$ embedded as $x\in\R^n \mapsto x+i0 \in \C^n=H$. Then two vectors $x,y\in V$ are 
$\sim$ equivalent if and only if $x=y$ or $x=-y$. Similarly, the projective space $\PP(V)$ is diffeomorphically equivalent to the real projective space $\R\PP^{n-1}$ which is of (real) dimension $n-1$. Thus
\[ \dim_{\R} (\Vo) = n. \]
\item Complex case. $V=\C^n$ which has real dimension $2n$. 
Then the projective space $\PP(V)=\C\PP^{n-1}$ has real dimension $2n-2$ 
(it is also a Kh\"{a}ler manifold) and thus
\[ \dim_{\R} (\Vo) = 2n-1. \]
\end{itemize}

The significance of the real dimension of $\Vo$ is encoded in the following result:
\begin{theorem}[\cite{BCE06}] 
If $m\geq 1+\dim_{\R}(\Vo)$ then for a (Zariski) generic frame $\fc$ of $m$ elements, the set of vectors $x\in V$ such that 
$\alpha^{-1}(\alpha(\hat{x}))$ has one point in $\hat{V}$ has dense interior in $V$.
\end{theorem}
The real case of this result is contained in Theorem 2.9, whereas the complex case is contained in Theorem 3.4. Both can be found in \cite{BCE06}.

\subsection{$\SS^{p,q}$} 

Consider now  
$\Sym(H)=\{T:\C^n\rightarrow \C^n~,~T=T^*\}$, the real vector space of self-adjoint operators over $H=\C^n$ endowed with the Hilbert-Schmidt scalar product
\index{product! Hilbert-Schmidt} $\ip{T}{S}_{HS}=\trace(TS)$. 
We also use the notation $\Sym(W)$ for the real vector space of symmetric operators over a (real or complex) vector space $W$. 
In both cases self-adjoint means the operator $T$ satisfies $\ip{Tx}{y}=\ip{x}{Ty}$ for every $x,y$ in the underlying vector space $W$.
$T^*$ means the adjoint operator of $T$, and therefore the transpose conjugate of $T$, when $T$ is
a matrix. When $T$ is a an operator acting on a real vector space, $T^T$ denotes its adjoint. For two vectors $x,y\in\C^n$ we denote 
\begin{equation}
\outp{x}{y}=\frac{1}{2}(xy^*+yx^*) \in \Sym(\C^n),
\end{equation}
their symmetric outer product\index{product! symmetric outer, \mbox{$\outp{\cdot}{\cdot}$}}.
On $\Sym(H)$ and $B(H)=\C^{n\times n}$ we consider the class of $p$-norms \index{norm! $p$-norm} defined by the
$p$-norm of the vector of singular values:
\begin{equation} \norm{T}_p = \left\{ 
\begin{array}{ccc}
\mbox{$\max_{1\leq k\leq n} \sigma_k(T)$} & for & \mbox{$p=\infty$} \\
\mbox{$\left( \sum_{k=1}^n \sigma_k^p\right)^{1/p} $} & for & \mbox{$1\leq p<\infty$}
\end{array}
\right.,
\end{equation}
where $\sigma_k=\sqrt{\lambda_k(T^*T)}$, ${1\leq k\leq n}$, are the singular values of $T$, with $\lambda_k(S)$, $1\leq k\leq n$, denoting the eigenvalues of $S$.

Fix two integers $p,q\geq 0$ and set
\begin{eqnarray}
\SS^{p,q}(H) & = & \{ T\in \Sym(H)~,~T~{\rm has~at~most}~ p~{\rm strictly~positive~ eigenvalues~} \nonumber \\
 & & {\rm and~ at~ most}~ q~{\rm strictly~negative~ eigenvalues} \},
\end{eqnarray}
\begin{eqnarray}
\mathring{\SS}^{p,q}(H) & = & \{ T\in \Sym(H)~,~T~{\rm has~exactly}~ p~{\rm strictly~positive~ eigenvalues~} \nonumber \\
& &  {\rm and~ exactly}~ q~{\rm strictly~negative~ eigenvalues} \}.
\end{eqnarray}
For instance $\mathring{\SS}^{0,0}(H)=\SS^{0,0}(H)=\{ 0\}$ and $\mathring{\SS}^{1,0}(H)$ is the set of all non-negative operators of rank exactly one. 
When there is no confusion we shall drop the underlying vector space $H=\C^n$ from notation. 

The following basic properties can be found in \cite{Bal13a}, Lemma 3.6; in fact, the last statement is a special instance of the Witt's decomposition theorem.
\begin{lemma}\mbox{}

\begin{enumerate}
\item For any $p_1\leq p_2$ and $q_1\leq q_2$, $\SS^{p_1,q_1}\subset \SS^{p_2,q_2}$;
\item For any nonnegative integers $p,q$ the following disjoint decomposition holds true
\begin{equation}
\SS^{p,q}=\cup_{r=0}^p\cup_{s=0}^q \mathring{\SS}^{r,s},
\end{equation}
where by convention $\mathring{\SS}^{p,q}=\emptyset$ for $p+q>n$.
\item For any $p,q\geq 0$,
\begin{equation}
-\SS^{p,q} = \SS^{q,p}.
\end{equation}
\item For any linear operator $T:H\rightarrow H$ (symmetric or not, invertible or not) and nonnegative integers $p,q$,
\begin{equation}
T^* \SS^{p,q}T \subset \SS^{p,q}.
\end{equation}
However if $T$ is invertible then $T^*\SS^{p,q}T=\SS^{p,q}$.
\item For any nonnegative integers $p,q,r,s$,
\begin{equation}
\SS^{p,q}+\SS^{r,s} = \SS^{p,q} - \SS^{s,r} = \SS^{p+r,q+s}.
\end{equation}
\end{enumerate}
\end{lemma}

The spaces $\SS^{1,0}$ and $\SS^{1,1}$ play a special role in the following section. Next we summarize
their properties (see Lemmas 3.7 and 3.9 in \cite{Bal13a}, and the comment after Lemma 9 in \cite{BCMN13a}). 

\begin{lemma}[{\rm Space $\SS^{1,0}$}] The following statements hold true:
\begin{enumerate}
\item $\mathring{\SS}^{1,0} =\{xx^*~,~x\in H,x\neq 0\}$;
\item $\SS^{1,0} =\{ xx^* ~,~x\in H\} = \{0\} \cup \{xx^* ~,~x\in H,x\neq 0\}$;
\item The set $\mathring{\SS}^{1,0}$ is a real analytic manifold in $\Sym(\C^n)$ of real dimension $2n-1$. 
As a real manifold, its tangent space at $X=xx^*$ is given by
\begin{equation}
T_X \mathring{\SS}^{1,0} = \left\{ \outp{x}{y} = \frac{1}{2}(xy^*+yx^*) ~,~y\in \C^n \right\}.
\end{equation}
The $\R$-linear embedding $\C^n \mapsto T_X \mathring{\SS}^{1,0}$ given by $y\mapsto \outp{x}{y}$ has
null space $\{ iax~,~a\in\R\}$.
\end{enumerate}
\end{lemma}

\begin{lemma}[{\rm Space $\SS^{1,1}$}] \label{l2.3}The following statements hold true:
\begin{enumerate}
\item $\SS^{1,1}=\SS^{1,0}-\SS^{1,0}=\SS^{1,0}+\SS^{0,1}=\{\outp{x}{y}~,~x,y\in H\}$;
\item For any vectors $x,y,u,v\in H$,
\begin{eqnarray}
 xx^*-yy^* & = & \outp{x+y}{x-y}= \outp{x-y}{x+y}, \\
\outp{u}{v} & = & \frac{1}{4} (u+v)(u+v)^* - \frac{1}{4}(u-v)(u-v)^*.
\end{eqnarray}
Additionally, for any $T\in\SS^{1,1}$ let $T=a_1 e_1e_1^* - a_2 e_2e_2^*$ be its spectral factorization with
$a_1,a_2\geq 0$ and $\ip{e_i}{e_j}=\delta_{i,j}$. Then
\[ T = \outp{\sqrt{a_1}e_1+\sqrt{a_2}e_2}{\sqrt{a_1}e_1-\sqrt{a_2}e_2}. \]
\item The set $\mathring{\SS}^{1,1}$ is a real analytic manifold in $\Sym(\C^n)$ of real dimension $4n-4$. Its tangent space at
 $X=\outp{x}{y}$ is given by
\begin{equation}
T_X \mathring{\SS}^{1,1} = \{ \outp{x}{u}+\outp{y}{v} = \frac{1}{2} (xu^*+ux^*+yv^*+vy^*)~,~u,v\in \C^n\}.
\end{equation}
The $\R$-linear embedding $\C^n\times \C^n\mapsto T_X \mathring{\SS}^{1,1}$ given by $(u,v)\mapsto \outp{x}{u}+
\outp{y}{v}$ has null space $\{ a(ix,0)+b(0,iy)+c(y,-x)+d(iy,ix)~,~a,b,c,d\in\R\}$. 
\item Let $T=\outp{u}{v}\in\SS^{1,1}$. Then its eigenvalues and $p$-norms are:
\begin{eqnarray}
a_{+} & = & \frac{1}{2} \left( \real(\ip{u}{v}) + \sqrt{\norm{u}^2\norm{v}^2 - (\imag(\ip{u}{v}))^2} \right) \geq 0, \\
a_{-} & = &  \frac{1}{2} \left( \real(\ip{u}{v}) - \sqrt{\norm{u}^2\norm{v}^2 - (\imag(\ip{u}{v}))^2} \right) \leq 0, \\
\norm{T}_1 & = & \sqrt{\norm{u}^2\norm{v}^2 - (\imag(\ip{u}{v}))^2}, \label{eq2.19} \\
\norm{T}_2 & = & \sqrt{ \frac{1}{2} \left(  \norm{u}^2\norm{v}^2 + (\real(\ip{u}{v}))^2 - (\imag(\ip{u}{v}))^2 \right)}, \\
\norm{T}_{\infty} & = &  \frac{1}{2} \left( |\real(\ip{u}{v})| + \sqrt{\norm{u}^2\norm{v}^2 - (\imag(\ip{u}{v}))^2} \right).
\end{eqnarray}
\item Let $T=xx^*-yy^*\in\SS^{1,1}$. Then its eigenvalues and $p$-norms are:
\begin{eqnarray}
a_{+} & = & \frac{1}{2} \left( \norm{x}^2 - \norm{y}^2 + \sqrt{(\norm{x}^2+\norm{y}^2)^2 - 4|\ip{x}{y}|^2} \right), \\
a_{-} & = & \frac{1}{2} \left( \norm{x}^2 - \norm{y}^2 - \sqrt{(\norm{x}^2+\norm{y}^2)^2 - 4|\ip{x}{y}|^2} \right), \\
\norm{T}_1 & = & \sqrt{(\norm{x}^2+\norm{y}^2)^2 - 4|\ip{x}{y}|^2}, \\
\norm{T}_2 & = & \sqrt{\norm{x}^4 + \norm{y}^4 - 2|\ip{x}{y}|^2}, \\
\norm{T}_{\infty} & = & \frac{1}{2} \left( |\norm{x}^2 - \norm{y}^2| + \sqrt{(\norm{x}^2+\norm{y}^2)^2 - 4|\ip{x}{y}|^2} \right).
\end{eqnarray}
\end{enumerate}
\end{lemma}

Note the above results hold true for the case of symmetric operators over the real subspace $V$.
In particular the factorization at Lemma \ref{l2.3}(1) implies that
\begin{equation}
\label{eq:2.27}
\SS^{1,1}(V) = \SS^{1,0}(V) - \SS^{1,0}(V) = \SS^{1,0}(V) + \SS^{0,1}(V) = \{ \outp{u}{v}~,~u,v\in V\}.
\end{equation}
More generally this result holds for subsets $V\subset H$ that are closed under addition and subtraction (such as modules over $\Z$). 

\subsection{Metrics}
The space $\hat{H}=\widehat{\C^n}$ admits two classes of distances (metrics). The first class is the ``natural metric" induced by the quotient space structure. The second metric is a matrix norm-induced distance.

Fix $1\leq p\leq \infty$.

The {\em natural metric} \index{distance! natural metric} denoted by $D_p:\hat{H}\times\hat{H}\rightarrow\R$ is defined by
\begin{equation}
D_p(\hat{x},\hat{y}) = \min_{\varphi \in [0,2\pi)} \norm{x-e^{i\varphi} y }_p,
\end{equation}
where $x\in\hat{x}$ and $y\in\hat{y}$. 
In the case $p=2$ the distance becomes
\[
D_2(\hat{x},\hat{y}) =  \sqrt{\norm{x}^2 + \norm{y}^2 - 2|\ip{x}{y}|}.
\]
By abuse of notation we use also $D_p(x,y)=D_p(\hat{x},\hat{y})$ since the distance does not depend on the choice of representatives.

The {\em matrix norm-induced distance} \index{distance! matrix-norm induced} denoted by $d_p:\hat{H}\times\hat{H}\rightarrow\R$ is defined by
\begin{equation}
d_p(\hat{x},\hat{y}) = \norm{xx^*-yy^*}_p,
\end{equation}
where again $x\in\hat{x}$ and $y\in\hat{y}$. In the case $p=2$ we obtain
\[ d_2(x,y) = \sqrt{\norm{x}^4 + \norm{y}^4 - 2|\ip{x}{y}|^2}. \]
By abuse of notation we use also $d_p(x,y)=d_p(\hat{x},\hat{y})$ since again the distance does not depend on the choice of representatives.

As analyzed in \cite{BZ14}, Proposition 2.4,  
$D_p$ is not Lipschitz equivalent to $d_p$, however $D_p$ is an equivalent distance to $D_q$
and similarily,  $d_p$ is equivalent to $d_q$, for any $1\leq p,q\leq q$ (see also \cite{BZ15} for the last claim below):

\begin{lemma}\mbox{}\label{l2.4}
\begin{enumerate}
\item For each $1\leq p\leq \infty$, $D_p$ and $d_p$ are distances (metrics) on $\hat{H}$;
\item $(D_p)_{1\leq p\leq \infty}$ are equivalent distances; that is, each $D_p$ induces the same topology on $\hat{H}$
and, for every $1\leq p,q\leq \infty$, the identity map $i:(\hat{H},D_p)\rightarrow (\hat{H},D_q)$, $i(x)=x$, is
Lipschitz continuous with Lipschitz constant
\[ \Lip^D_{p,q,n} = \max(1,n^{\frac{1}{q}-\frac{1}{p}}). \]
\item $(d_p)_{1\leq p\leq \infty}$ are equivalent distances, that is, each $d_p$ induces the same topology on $\hat{H}$
and, for every $1\leq p,q\leq \infty$, the identity map $i:(\hat{H},d_p)\rightarrow (\hat{H},d_q)$, $i(x)=x$, is
Lipschitz continuous with Lipschitz constant
\[ \Lip^d_{p,q,n} = \max(1,2^{\frac{1}{q}-\frac{1}{p}}). \]
\item The identity map $i:(\hat{H},D_p)\rightarrow (\hat{H},d_p)$, $i(x)=x$ is continuous, but it is not Lipschitz continuous. The identity map $i:(\hat{H},d_p)\rightarrow (\hat{H},D_p)$, $i(x)=x$ is continuous but it is not Lipschitz continuous. Hence the induced topologies on $(\hat{H},D_p)$ and $(\hat{H},d_p)$ are the same, but the 
corresponding distances are not Lipschitz equivalent.
\item The metric space $(\hat{H},d_p)$ is isometrically isomorphic to $\SS^{1,0}$ endowed with the $p$-norm.
The isomorphism is given by the map
\[ \kappa_\beta:\hat{H}\rightarrow \SS^{1,0}~~,~~x\mapsto \outp{x}{x}=xx^*. \] 
\item The metric space $(\hat{H},D_2)$ is Lipschitz isomorphic (not isometric) with $\SS^{1,0}$ endowed with the $2$-norm. The bi-Lipschitz map
\[ \kappa_{\alpha}:\hat{H}\rightarrow \SS^{1,0}~~,~~x\mapsto \kappa_{\alpha}(x)=\left\{ 
\begin{array}{rcl}
\mbox{$\frac{1}{\norm{x}}xx^*$} & if & \mbox{$x\neq 0$} \\
0 & & otherwise
\end{array} \right. 
\]
has lower Lipschitz constant $1$ and upper Lipschitz constant $\sqrt{2}$.
\end{enumerate}
\end{lemma}

Note the Lipschitz constant $\Lip^D_{p,q,n}$ is equal to the operator norm of the identity map between $(\C^n,\norm{\cdot}_p)$ and $(\C^n,\norm{\cdot}_q)$:
 $\Lip^D_{p,q,n}=\norm{I}_{l^p(\C^n)\rightarrow l^q(\C^n)}$.
Note also the equality $\Lip^d_{p,q,n}=\Lip^D_{p,q,2}$.
A consequence of the last two claims in the above result is that while the identity map between $(\hat{H},D_p)$ and $(\hat{H},d_q)$  is not 
bi-Lipschitz, the map $x\mapsto \frac{1}{\sqrt{\norm{x}}}x$ is bi-Lipschitz.

\section{The Injectivity Problem}

In this section we summarize existing results on the injectivity of the maps $\alpha$ and $\beta$.
Our plan is to present the real and the complex case in a unified way.

Recall that $V$ is a real vector space which is also a subset of $H=\C^n$. 
The special two cases are $V=\R^n$ (the real case) and $V=\C^n$ (the complex case). 

First we describe the {\em realification procedure} \index{realification} of $H$ and $V$. Consider the $\R$-linear map $\j:\C^n\rightarrow\R^{2n}$
defined by
\[ 
\j(x) = \left[ \begin{array}{c}
\mbox{$\real(x)$} \\
\mbox{$\imag(x)$}  \end{array} \right].
\]
Let $\Vc = \j(V)$ be the embedding of $V$ into $\R^{2n}$, and let $\Pi$ denote the orthogonal projection (with respect to the real scalar product on $\R^{2n}$) onto $\Vc$.
Let $J$ denote the folowing orthogonal antisymmetric $2n\times 2n$ matrix
\begin{equation}
J = \left[ \begin{array}{cc}
\mbox{$0$} & \mbox{$-I_n$} \\
\mbox{$I_n$} & \mbox{$0$} 
\end{array} \right] ,
\end{equation}
where $I_n$ denotes the $n\times n$ identity matrix.
Note that $J^T=-J$, $J^2=-I_{2n}$ and $J^{-1}=-J$.

Each vector $f_k$ of the frame set $\fc=\{f_1,\ldots,f_m\}$ is mapped by $\j$ onto a vector in $\R^{2n}$ denoted by $\varphi_k$,
 and a symmetric operator in $\SS^{2,0}(\R^{2n})$ denoted by $\Phi_k$:
\begin{equation}
\varphi_k = \j(f_k) = \left[ 
\begin{array}{c}
\mbox{$\real(f_k)$} \\ \mbox{$\imag(f_k)$} 
\end{array} \right] ~~,~~ \Phi_k = \varphi_k \varphi_k^T + J\varphi_k \varphi_k^TJ^T.
\end{equation}
Note that when $f_k\neq 0$ the symmetric form $\Phi_k$ has rank 2  and belongs to $\mathring{\SS}^{2,0}$.
Its spectrum has two distinct eigenvalues: $\norm{\varphi_k}^2=\norm{f_k}^2$ with multiplicity $2$,
and $0$ with multiplicity $2n-2$. Furthermore, $\frac{1}{\norm{\varphi_k}^2} \Phi_k$ is a rank 2 projection.

Let $\xi=\j(x)$ and $\eta=\j(y)$ denote the realifications of vectors $x,y\in\C^n$. Then a bit of algebra shows that
\begin{eqnarray}
\ip{x}{f_k} & = & \ip{\xi}{\varphi_k} + i \ip{\xi}{J\varphi_k}, \\
\nonumber
\ip{F_k}{xx^*}_{HS}=\trace\left(F_k xx^* \right) = |\ip{x}{f_k}|^2 = \ip{\Phi_k \xi}{\xi} & = & \ip{\Phi_k}{\xi\xi^T}_{HS},\\
\ip{F_k}{\outp{x}{y}}_{HS}=\trace\left( F_k \outp{x}{y} \right) = \real(\ip{x}{f_k}\ip{f_k}{y}) 
& = &\ip{\Phi_k \xi}{\eta} = \trace(\Phi_k \outp{\xi}{\eta}) \nonumber \\
& =& \ip{\Phi_k}{\outp{\xi}{\eta}}_{HS},
\nonumber
\end{eqnarray}
where $F_k = \outp{f_k}{f_k} = f_kf_k^*\in\SS^{1,0}(H)$.

The following objects play an important role in the subsequent theory:
\begin{eqnarray}
& & \Rbc:\C^n \rightarrow \Sym(\C^n) ~,~ \Rbc(x) = \sum_{k=1}^m |\ip{x}{f_k}|^2 f_kf_k^* ~,~x\in\C^n,\\
& & \Rsc:\R^{2n}\rightarrow \Sym(\R^{2n}) ~,~ \Rsc(\xi) = \sum_{k=1}^m \Phi_k \xi \xi^T \Phi_k ~,~\xi\in\R^{2n},\\
& & \Ssc:\R^{2n}\rightarrow \Sym(\R^{2n}) ~,~ \Ssc(\xi) = \sum_{k:\Phi_k\xi\neq 0}\frac{1}{\ip{\Phi_k\xi}{\xi}}
\Phi_k\xi\xi^T\Phi_k~,~\xi\in\R^{2n}, \\
& & \Zc:\R^{2n} \rightarrow\R^{2n\times m} ~,~ \Zc(\xi) = \left[ 
\begin{array}{ccccc} \mbox{$\Phi_1\xi$} & \mbox{$|$} & \mbox{$\ldots$} & \mbox{$|$} & \mbox{$\Phi_m \xi$} \end{array} \right] ~,~\xi\in\R^{2n}.
\end{eqnarray}
Note $\Rsc = \Zc \Zc^T$.

 Following \cite{BBCE07} we note that $|\ip{x}{f_k}|^2$ is the Hilbert-Schmidt scalar product between two rank 1
symmetric forms: 
\[ |\ip{x}{f_k}|^2  = \trace \left(F_k X \right) =  \ip{F_k}{X}_{HS},  \]
where $X=xx^*$. Thus the nonlinear map $\beta$ induces a linear map on the real vector space $\Sym(\C^n)$ of
symmetric forms over $\C^n$:
\begin{equation}\label{eq:Ac}
\Ac: \Sym(\C^n) \rightarrow \R^m ~~,~~(\Ac(T))_k = \ip{T}{F_k}_{HS} = \ip{T f_k}{f_k} ~,~{1\leq k\leq m}
\end{equation}
Similarly it induces a linear map on $\Sym(\R^{2n})$, the space of symmetric forms over $\R^{2n}=\j(\C^n)$,
that is denoted by $\Acr$:
\begin{equation}
\Acr: \Sym(\R^{2n})\rightarrow\R^{m} ~~,~~(\Acr(T))_k = \ip{T}{\Phi_k}_{HS}= 
\ip{T\varphi_k}{\varphi_k}+\ip{TJ\varphi_k}{J\varphi_k}~,~{1\leq k\leq m}.
\end{equation}
Now we are ready to state a necessary and sufficient condition for injectivity that works in both the real and the complex case:

\begin{theorem}[\cite{HMW11,BCMN13a,Bal13a}] Let $H=\C^n$ and let $V$ be a real vector space that is also a subset of $H$, $V\subset H$. 
Denote by $\Vc=\j(V)$ the realification of $V$. Assume $\fc$ is a frame for $V$. The following statements are equivalent:
\begin{enumerate}
\item The frame $\fc$ is phase retrievable with respect to $V$;
\item $\ker\Ac \cap \left(\Sonezero(V)-\Sonezero(V) \right) =\{0\}$;
\item $\ker\Ac \cap \Soneone(V) =\{0\}$;
\item $\ker\Ac\cap (\SS^{2,0}(V)\cup\Soneone(V)\cup\SS^{0,2}) = \{0\} $;
\item There do not exist vectors $u,v\in V$ with $\outp{u}{v}\neq 0$ so that
\[ \real\left( \ip{u}{f_k}\ip{f_k}{v} \right) = 0~~,~~\forall\, 1\leq k\leq m; \]
\item $\ker\Acr \cap \left(\Sonezero(\Vc)-\Sonezero(\Vc) \right) =\{0\}$;
\item $\ker\Acr \cap \Soneone(\Vc) =\{0\}$;
\item There do not exist vectors $\xi,\eta\in \Vc$, with $\outp{\xi}{\eta}\neq 0$ so that
\[ \ip{\Phi_k\xi}{\eta} = 0~~,~~\forall\, 1\leq k\leq m. \]
\end{enumerate}\label{t3.1}
\end{theorem}

{\bf Proof}. 

$(1)\Leftrightarrow (2)$ It is immediate once we notice that any element in the null space of $\Ac$ of the 
form $xx^*-yy^*$ 
implies $\Ac(xx^*)=\Ac(yy^*)$ for some $x,y\in V$ with $\hat{x}\neq \hat{y}$. 

$(2)\Leftrightarrow(3)$ and $(3)\Leftrightarrow(5)$ are consequences of (\ref{eq:2.27}).

For $(4)$ first note that $\ker\Ac\cap \SS^{2,0}(V)=\{0\}=\ker\Ac\cap\SS^{0,2}(V)$ since $\fc$ is a frame for $V$. Thus $(3)\Leftrightarrow(4)$.

(6),(7) and (8) are simply restatements of (2),(3) and (4) using the realification procedure.
$\Box$

In the case $(4)$ above, note that $\SS^{2,0}(V)\cup\Soneone(V)\cup\SS^{0,2}(V)$ is the set of all rank less than or equal to 2 symmetric operators in $\Sym(V)$ (This statement has
been proposed in \cite{BCMN13a}).

The above general injectivity result is next made more explicit in the cases $V=\C^n$ and $V=\R^n$.

\begin{theorem}[\cite{BCE06,Bal12a}](The real case) Assume $\fc\subset\R^{n}$. The following are equivalent:
\begin{enumerate}
\item $\fc$ is phase retrievable for $V=\R^n$;
\item $\Rbc(x)$ is invertible for every $x\in\R^n$, $x\neq 0$;
\item There do not exist vectors $u,v\in\R^n$ with $u\neq 0$ and $v\neq 0$ so that
\[ \ip{u}{f_k}\ip{f_k}{v} = 0~~,~~\forall\, 1\leq k\leq m; \]
\item For any disjoint partition of the frame set $\fc=\fc_1\cup\fc_2$, either $\fc_1$ spans $\R^n$ or $\fc_2$ spans $\R^n$.
\end{enumerate}
\label{th3.2}
\end{theorem}

Recall a set $\fc\subset \C^n$ is called {\em full spark} \index{full spark} if any subset of $n$ vectors is linearly independent.
Then an immediate corollary of the above result is the following
\begin{corollary}[\cite{BCE06}] Assume $\fc\subset\R^n$. Then
\begin{enumerate}
\item If $\fc$ is phase retrievable for $\R^n$ then $m\geq 2n-1$;
\item If $m=2n-1$, then $\fc$ is phase retrievable if and only if $\fc$ is full spark;
\end{enumerate}
\end{corollary}
{\bf Proof}. 

Indeed, the first claim follows from Theorem \ref{th3.2}{(4): If $m\leq 2n-2$ then there is a partition of $\fc$ into two subsets each of cardinal less than or equal to $n-1$. Thus neither set can span $\R^n$. Contradiction.

The second claim is immediate from the same statement as above. $\Box$

A more careful analysis of Theorem \ref{th3.2}(4) gives a recipe of constructing two non-similar vectors $x,y\in\R^{n}$ so that $\alpha(x)=\alpha(y)$. Indeed, if $\fc=\fc_1\cup\fc_2$ so that $\dim{\span}(\fc_1)<n$ and $\dim{\span}(\fc_2)<n$ then there are non-zero vectors $u,v\in\R^n$ with $\ip{u}{f_k}=0$ for all $k\in I$, and $\ip{v}{f_k}=0$ for all $k\in I^c$. Here $I$ is the index set of frame vectors in $\fc_1$ and $I^c$ denotes its complement in $\{1,\ldots,m\}$. Set $x=u+v$ and $y=u-v$. Then $|\ip{x}{f_k}|=|\ip{v}{f_k}|=|\ip{y}{f_k}|$ for all $k\in I$, and $|\ip{x}{f_k}|=|\ip{u}{f_k}|=|\ip{y}{f_k}|$ for all $k\in I^c$. Thus $\alpha(x)=\alpha(y)$, but $x\neq y$ and $x\neq -y$.

\begin{theorem}[\cite{BCMN13a,Bal13a}](The complex case) The following are equivalent:
\begin{enumerate}
\item $\fc$ is phase retrievable for $H=\C^n$;
\item $\rank(\Zc(\xi))=2n-1$ for all $\xi\in\R^{2n}$, $\xi\neq 0$;
\item $\dim\ker \Rsc(\xi) = 1$ for all $\xi\in\R^{2n}$, $\xi\neq 0$;
\item There do not exist $\xi,\eta \in\R^{2n}$, $\xi\neq 0$ and $\eta\neq 0$ so that $\ip{J\xi}{\eta}= 0$ and
\begin{equation}
\ip{\Phi_k\xi}{\eta} = 0~~,~~\forall 1\leq k\leq m.
\end{equation}
\end{enumerate}
\end{theorem}

In terms of cardinality, here is what we know:
\begin{theorem}[\cite{Mil67,HMW11,BH13,Bal13b,MV13,CEHV13,Viz15}] \mbox{}
\begin{enumerate}
\item \cite{HMW11} If $\fc$ is a phase retrievable frame for $\C^n$ then 
\begin{equation}
\label{eq:HMW11}
m\geq 4n-2-2b + \left\{ 
\begin{array}{rl}
\mbox{$2$} & \mbox{if $n$ odd and $b=3\,mod\,4$} \\
\mbox{$1$} & \mbox{if $n$ odd and $b=2\,mod\,4$} \\
\mbox{$0$} & \mbox{otherwise}
\end{array} \right.,
\end{equation}
where $b=b(n)$ denotes the number of 1's in the binary expansion of $n-1$.
\item \cite{BH13} For any positive integer $n$ there is a frame with $m=4n-4$ vectors so that $\fc$ is phase retrievable for $\C^n$;
\item \cite{CEHV13} If $m\geq 4n-4$ then a (Zariski) generic frame is phase retrievable for $\C^n$;
\item \cite{Bal13b} The set of phase retrievable frames is open in $\C^{n}\times\cdots\times\C^n$. In particular phase retrievable property is stable under small perturbations.   
\item \cite{CEHV13} If $n=2^k+1$ and $m\leq 4n-5$ then $\fc$ cannot be phase retrievable for $\C^n$.
\item \cite{Viz15} For $n=4$ there is a frame with $m=11<4n-4=12$ vectors that is phase retrievable for $\C^n$.
\end{enumerate}
\end{theorem}

\section{Robustness of Reconstruction}

In this section we analyze stability bounds for reconstruction. Specifically we analyze two types of margins:
\begin{itemize}
\item Deterministic, worst-case type bounds: These bounds are given by lower Lipschitz constants of the forward nonliner analysis maps;
\item Stochastic, average type bounds: Cramer-Rao Lower Bounds (CRLB).
\end{itemize}

\subsection{Bi-Lipschitzianity of the Nonlinear Analysis Maps}

In Section 2 we introduced two distances on $\hat{H}$. As the following theorem shows, the nonlinear maps
$\alpha$ and $\beta$ are bi-Lipschitz \index{Lipschitz! bi-Lipschitz} with respect to the corresponding distance:

\begin{theorem}\cite{Bal12a,EM12,BCMN13a,Bal13a,BW13,BZ14,BZ15a,BZ15}\label{t4.1}
Let $\fc$ be a phase retrievable frame for $V$, a real linear space, subset of $H=\C^n$. Then:
\begin{enumerate}
\item The nonlinear map $\alpha:(\hat{V},D_2)\rightarrow (\R^m,\norm{}_2)$ is bi-Lipschitz. Specifically there are positive constants $0<A_0\leq B_0<\infty$ so that
\begin{equation}\label{eq:alpha}
\sqrt{A_0} D_2(x,y) \leq \norm{\alpha(x)-\alpha(y)}_2 \leq \sqrt{B_0} D_2(x,y) ~~,~~\forall\, x,y \in V.
\end{equation}
\item The nonlinear map $\beta:(\hat{V},d_1)\rightarrow (\R^m,\norm{}_2)$ is bi-Lipschitz. Specifically there are positive constants $0<a_0\leq b_0<\infty$ so that
\begin{equation}\label{eq:beta}
\sqrt{a_0} d_1(x,y) \leq \norm{\beta(x)-\beta(y)}_2 \leq \sqrt{b_0} d_1(x,y) ~~,~~\forall\, x,y \in V.
\end{equation} 
\end{enumerate}
The converse is also true: If either (\ref{eq:alpha}) or (\ref{eq:beta}) holds true for all $x,y\in V$ then $\fc$ is phase retrievable for $V$.
\end{theorem}
The choice of distance $D_2$ and $d_1$ in the statement of this theorem is only for reasons of convenience since these specific constants will appear later in the text. 
Any other distance $D_p$ instead of $D_2$, and $d_q$ instead of $d_1$ would work. 
The Lipschitz constants would be different, of course. This result was first obtained for the real case in \cite{EM12} for the map $\alpha$ and in \cite{Bal12a} for the map $\beta$. The complex case for map $\beta$ was shown independently in \cite{BCMN13a} and \cite{Bal13a}. The complex case for the more challenging map $\alpha$ was proved in \cite{BZ15}. The paper \cite{BW13} computes the optimal bound $A_0$ in the real case. The statement presented here (Theorem \ref{t4.1}) unifies these two cases.

On the other hand the condition that $\fc$ is phase retrievable for $V$ is equivalent to the existence of
a lower bound for a family of quadratic forms.
We state this condition now:
\begin{theorem}\label{t4.2}
Let $\fc\subset H=\C^n$ and let $V$ be a real vector space, subset of $H$. Denote by $\Vc=\j(V)\subset\R^{2n}$
 the realification of $V$, and let $\Pi$ denote the projection onto $\Vc$. Then the following statements are equivalent:
\begin{enumerate}
\item $\fc$ is phase retrievable for $V$;
\item There is a constant $a_0>0$ so that
\begin{equation}
\label{eq:a0-1}
 \Pi \Rsc(\xi) \Pi \geq a_0 \Pi P_{J\xi}^{\perp} \Pi ~~,~~\forall \, \xi\in\Vc,\norm{\xi}=1, 
\end{equation}
where $P_{J\xi}^{\perp} = I_{2n}-P_{J\xi}=I_{2n}-J\xi\xi^TJ^T$ is the orthogonal projection onto the orthogonal complement to $J\xi$;
\item There is $a_0>0$ so that for all $\xi,\eta\in\R^{2n}$,
\begin{equation}
\label{eq:a0-2}
\sum_{k=1}^m |\ip{\Pi\Phi_k\Pi\xi}{\eta}|^2 \geq a_0 \left(\norm{\Pi\xi}^2\norm{\Pi\eta}^2-|\ip{J\Pi\xi}{\Pi\eta}|^2 \right).
\end{equation}
\end{enumerate}
\end{theorem}
Note the same constant $a_0$ can be chosen in (\ref{eq:beta}) and (\ref{eq:a0-1}) and (\ref{eq:a0-2}).
This result was shown separately for the real and complex case. Here we state these conditions in a unified way.

{\bf Proof}.

$(1)\Leftrightarrow(2)$ If  $\fc$ is a phase retrievable frame for $V$ then, by Theorem \ref{t3.1}(8), 
for all vectors $\xi,\eta\in \Vc$, with $\outp{\xi}{\eta}\neq 0$ we have $\ip{\Phi_k\xi}{\eta} \neq 0$, for some $1\leq k\leq m$. Take $\mu\in\R^{2n}$ and set $\eta=\Pi\mu$. Normalize $\xi$ to
$\norm{\xi}=1$. Then
\[ \sum_{k=1}^m |\ip{\Phi_k\xi}{\eta}|^2 = \ip{\Rsc(\xi)\Pi\mu}{\Pi\mu}, \]
and by (\ref{eq2.19}),
\[ \norm{\outp{\xi}{\eta}}_1^2 = \norm{\xi}^2\norm{\eta}^2-|\ip{\xi}{J\eta}|^2 = \norm{\Pi\mu}^2 - |\ip{J\xi}{\Pi\mu}|^2 = \ip{(I_{2n}-J\xi\xi^TJ^T)\Pi\mu}{\Pi\mu}. \]
Thus if $\mu$ satisfies $\outp{\xi}{\Pi\mu}=0$ then it must also satisfy $\Pi\mu=tJ\xi$ 
for some real $t$. In this case $\Pi\mu$ lies in the null space of $\Rsc(\xi)$. 
In particular this proves that the following quotient of quadratic forms
\[ \frac{\ip{\Pi\Rsc(\xi)\Pi\mu}{\mu}}{\ip{\Pi(I_{2n}-J\xi\xi^TJ^T)\Pi\mu}{\mu}} \]
is bounded above and below away from zero. 
This proves that (\ref{eq:a0-1}) must hold for some $a_0>0$. 
Conversely, if (\ref{eq:a0-1}) holds true, then for every $\xi,\eta\in\Vc$ with $\outp{\xi}{\eta}\neq 0$, $\ip{P_{J\xi}^{\perp}\eta}{\eta}\neq 0$ and thus  $\ip{\Phi_k\xi}{\eta}\neq 0$ for some $k$.
This shows that $\fc$ is a phase retrievable frame for $V$.

$(2)\Leftrightarrow(3)$ This follows by writing out (\ref{eq:a0-1}) explicitly. 

$\Box$

\begin{remark}
Condition (2) of this theorem expressed by Equation (\ref{eq:a0-1}) can be used to check if a given frame is phase retrievable as we explain next.

In the real case, $\Pi=I_{n}\oplus 0$, and this condition reduces to
\[ \Rbc(x) = \sum_{k=1}^m |\ip{x}{f_k}|^2 f_kf_k^T \geq a_0 \norm{x}^2 I_{H}~~,~~\forall x\in H=\R^n. \]
In turn this is equivalent to any of the conditions of Theorem \ref{th3.2}.

In the complex case the condition (\ref{eq:a0-1}) turns into
\begin{equation}
\label{eq:a0-11}
 \lambda_{2n-1}(\Rsc (\xi))\geq a_0 ~~,~~\forall \xi\in\R^{2n},\norm{\xi}=1, 
\end{equation}
where $\lambda_{2n-1}(\Rsc(\xi))$ denotes the next to the smallest eigenvalue of $\Rsc(\xi)$.
The algorithm requires an upper bound for $b_0=\max_{\norm{\xi}=1}\lambda_1(\Rsc(\xi))$. For instance
$b_0\leq B\max_k\norm{f_k}^2$, where $B$ is the frame upper bound \cite{Bal13b}.
The condition (\ref{eq:a0-11}) can be checked using an $\eps$-net of the unit sphere in $\R^{2n}$. 
Specifically let $\{\xi_j^\eps\}$ be such an $\eps$-net, that is $\norm{\xi_j^{\eps}}=1$ and $\norm{\xi_j^{\eps}-\xi_k^{\eps}}<\eps$ for all $j\neq k$. 
Set  $a_0=\frac{1}{2}\min_j\lambda_{2n-1}(\Rsc(\xi_j^\eps))$.
If $2b_0\eps\leq a_0$ then stop, otherwise set $\eps=\frac{1}{2}\eps$ and construct a new $\eps$-net.

The condition $2 b_0\eps\leq a_0$ guarantees that for every $\xi\in\R^{2n}$ with $\norm{\xi}=1$, 
$\lambda_{2n-1}(\Rsc(\xi))\geq a_0$ since (see also \cite{Bal13b} for a similar derivation)
$$\norm{\Rsc(\xi)-\Rsc(\xi_j^\eps)}\leq \sqrt{b_0^2\norm{\xi-\xi^\eps_j}\norm{\xi+\xi_j^\eps}}\leq 2b_0\norm{\xi-\xi_j^\eps}\leq 2b_0 \eps$$
and by Weyl's perturbation theorem (see III.2.6 in \cite{Bhatia}) 
$$\lambda_{2n-1}(\Rsc(x))\geq \lambda_{2n-1}(\Rsc(\xi_j^\eps))-\norm{\Rsc(\xi)-\Rsc(\xi_j^\eps)}\geq 2a_0-2b\eps\geq a_0.$$

Unfortunately such an approach has at least an NP computational cost since the cardinality of an 
$\eps$-net is of the order $\left(\frac{1}{\eps}\right)^n$. 

\end{remark}

The computations of lower bounds is fairly subtle. In fact there is a distinction between local bounds and global bounds. Specifically for every $z\in V$ we define the following bounds:

The {\em type I local lower Lipschitz bounds} \index{Lipschitz! bounds} are defined by:
\begin{eqnarray}
A(z) & = & \lim_{r\rightarrow 0} \inf_{x,y\in V, D_2(x,z)<r,D_2(y,z)<r} \frac{\norm{\alpha(x)-\alpha(y)}_2^2}{D_2(x,y)^2}, \\
a(z) & = & \lim_{r\rightarrow 0} \inf_{x,y\in V, d_1(x,z)<r,d_1(y,z)<r} \frac{\norm{\beta(x)-\beta(y)}_2^2}{d_1(x,y)^2}. 
\end{eqnarray}
The {\em type II local lower Lipschitz bounds} are defined by:
\begin{eqnarray}
\tilde{A}(z) & = & \lim_{r\rightarrow 0} \inf_{y\in V, D_2(y,z)<r} \frac{\norm{\alpha(z)-\alpha(y)}_2^2}{D_2(z,y)^2}, \\
\tilde{a}(z) & = & \lim_{r\rightarrow 0} \inf_{y\in V, d_1(y,z)<r} \frac{\norm{\beta(z)-\beta(y)}_2^2}{d_1(z,y)^2}. 
\end{eqnarray}
Similarly the {\em type I local upper Lipschitz bounds} are defined by:
\begin{eqnarray}
B(z) & = & \lim_{r\rightarrow 0} \sup_{x,y\in V, D_2(x,z)<r,D_2(y,z)<r} \frac{\norm{\alpha(x)-\alpha(y)}_2^2}{D_2(x,y)^2}, \\
b(z) & = & \lim_{r\rightarrow 0} \sup_{x,y\in V, d_1(x,z)<r,d_1(y,z)<r} \frac{\norm{\beta(x)-\beta(y)}_2^2}{d_1(x,y)^2} 
\end{eqnarray}
and the {\em type II local upper Lipschitz bounds} are defined by:
\begin{eqnarray}
\tilde{B}(z) & = & \lim_{r\rightarrow 0} \sup_{y\in V, D_2(y,z)<r} \frac{\norm{\alpha(z)-\alpha(y)}_2^2}{D_2(z,y)^2}, \\
\tilde{b}(z) & = & \lim_{r\rightarrow 0} \sup_{y\in V, d_1(y,z)<r} \frac{\norm{\beta(z)-\beta(y)}_2^2}{d_1(z,y)^2}. 
\end{eqnarray}

The {\em global lower bounds} are defined by
\begin{eqnarray}
A_0& = & \inf_{x,y\in V, D_2(x,y)>0} \frac{\norm{\alpha(x)-\alpha(y)}_2^2}{D_2(x,y)^2}, \\
a_0 & = & \inf_{x,y\in V, d_1(x,y)>0} \frac{\norm{\beta(x)-\beta(y)}_2^2}{d_1(x,y)^2}, 
\end{eqnarray}
whereas the {\em global upper bounds} are defined by
\begin{eqnarray}
B_0 & = & \sup_{x,y\in V, D_2(x,y)>0} \frac{\norm{\alpha(x)-\alpha(y)}_2^2}{D_2(x,y)^2}, \\
b_0 & = & \sup_{x,y\in V, d_1(x,y)>0} \frac{\norm{\beta(x)-\beta(y)}_2^2}{d_1(x,y)^2} 
\end{eqnarray}
 and represent the square of the corresponding Lipschitz constants \index{Lipschitz! constants}.

Due to homogeneity $A_0=A(0)$, $B_0=B(0)$, $a_0=a(0)$, and $b_0=b(0)$. 
On the other hand, for $z\neq 0$, 
$A(z)=A(\frac{z}{\norm{z}})$, $B(z)=B(\frac{z}{\norm{z}})$, $a(z)=a(\frac{z}{\norm{z}})$, and 
$b(z)=b(\frac{z}{\norm{z}})$.
Note that $A(z)$ stands for the local lower Lipschitz bound of type I at $z$, whereas $A$ denotes the optimal lower frame bound of $\fc$. 

The exact expressions of these bounds are summarized by the following results. 
For any $I\subset \{1,2,\ldots,m\}$ let $\fc[I]=\{f_k~,~k\in I\}$ denote the subset indexed by $I$. 
Also let $\sigma_1^2[I]$ and  $\sigma_n^2[I]$ denote the upper and the lower frame bound of set $\fc[I]$, respectively. Thus:
\[ \sigma_1^2[I] = \lambda_{\textnormal{max}} \left(\sum_{k\in I}f_kf_k^*\right) ~,~ \sigma_n^2[I] = \lambda_{\textnormal{min}} \left(\sum_{k\in I}f_kf_k^*\right). \]
As usual, $I^c$ denotes the complement of the index set $I$, that is $I^c=\{1,\ldots,m\}\setminus I$. 

\begin{theorem}[\cite{BW13,BCMN13a}](The real case) Assume $\fc\subset\R^n$ is a phase retrievable frame for $\R^n$. 
Let $A$ and $B$ denote its optimal lower and upper frame bound, respectively. Then:
\begin{enumerate}
\item For every $0\neq x\in\R^n$, $A(x) = \sigma_n^2[\supp(\alpha(x))]$, where $\supp(\alpha(x))=\{k~,~\ip{x}{f_k}\neq 0\}$;
\item For every $x\in\R^n$, $\tilde{A}(x)=A$;
\item $A_0=A(0) = \min_{I}(\sigma_n^2[I]+\sigma_n^2[I^c]) $;
\item For every $x\in\R^n$, $B(x)=\tilde{B}(x)=B$;
\item $B_0=B(0) = \tilde{B}(0)=B$, the optimal upper frame bound;
\item For every $0\neq x\in\R^n$, $a(x) =\tilde{a}(x)=\lambda_{\textnormal{min}}(\Rbc(x))/\norm{x}^2$;
\item $a_0=a(0)=\tilde{a}(0)=\min_{\norm{x}=1} \lambda_{\textnormal{min}}(\Rbc(x))$;
\item For every $0\neq x\in\R^n$, $b(x)=\tilde{b}(x)=\lambda_{\textnormal{max}}(\Rbc(x))/\norm{x}^2$;
\item $b_0=b(0)=\tilde{b}(0)=\max_{\norm{x}=1}\lambda_{\textnormal{max}}(\Rbc(x))$;
\item $a_0$ is the largest constant so that
\[ \Rbc(x) \geq a_0 \norm{x}^2 I_n ~~,~~\forall\, x\in \R^n, \]
or, equivalently,
\[ \sum_{k=1}^m |\ip{x}{f_k}|^2 |\ip{y}{f_k}|^2 \geq a_0 \norm{x}^2 \norm{y}^2 ~~,~~\forall\, x,y\in\R^n; \]
\item $b_0$ is the $4^{th}$ power of the frame analysis operator norm $T:(\R^n,\norm{\cdot}_2)\rightarrow(\R^m,\norm{\cdot}_4)$,
\[ b_0 = \norm{T}_{B(l^2,l^4)}^4 = \max_{\norm{x}_2=1} \sum_{k=1}^m |\ip{x}{f_k}|^4. \]
\end{enumerate}
\end{theorem}


The complex case is subtler. The following result presents some of the local and global Lipschitz bounds.
\begin{theorem}[\cite{BZ15}](The complex case) Assume $\fc$ is phase retrievable for $H=\C^n$ and $A,B$ are its optimal frame bounds. 
Then:
\begin{enumerate}
\item For every $0\neq z\in\C^n$, $A(z)=\lambda_{2n-1}\left( \Ssc(\j(z)) \right)$ (the next to the smallest eigenvalue);
\item $A_0=A(0)>0$;
\item For every $z\in\C^n$, $\tilde{A}(z)=\lambda_{2n-1}\left( \Ssc(\j(z))+\sum_{k:\ip{z}{f_k}=0}\Phi_k \right)$ (the next to the smallest eigenvalue);
\item $\tilde{A}(0)=A$, the optimal lower frame bound;
\item For every $z\in\C^n$, $B(z) = \tilde{B}(z)=\lambda_1\left( \Ssc(\j(z))+\sum_{k:\ip{z}{f_k}=0}\Phi_k \right)$ (the largest eigenvalue);
\item $B_0=B(0)=\tilde{B}(0)=B$, the optimal upper frame bound;
\item For every $0\neq z\in\C^n$, $a(z)=\tilde{a}(z) = \lambda_{2n-1}(\Rsc(\j(z)))/\norm{z}^2$ (the next to the smallest eigenvalue);
\item For every $0\neq z\in\C^n$, $b(z) = \tilde{b}(z) = \lambda_{1}(\Rsc(\j(z)))/\norm{z}^2$ (the largest eigenvalue);
\item $a_0$ is the largest constant so that
\[ \Rsc(\xi) \geq a_0 (I-J\xi\xi^TJ^T)~~,~~\forall\, \xi\in\R^{2n},\norm{\xi}=1, \]
or, equivalently,
\[ \sum_{k=1}^m |\ip{\Phi_k\xi}{\eta}|^2 \geq a_0 \left( \norm{\xi}^2\norm{\eta}^2-|\ip{J\xi}{\eta}|^2 \right)
~~,~~\forall\, \xi,\eta\in\R^{2n}; \]
\item $b(0)=\tilde{b}(0)=b_0$ is the $4^{th}$ power of the frame analysis operator norm $T:(\C^n,\norm{\cdot}_2)\rightarrow(\R^m,\norm{\cdot}_4)$,
\[ b_0 = \norm{T}_{B(l^2,l^4)}^4 = \max_{\norm{x}_2=1} \sum_{k=1}^m |\ip{x}{f_k}|^4; \]
\item $\tilde{a}(0)$ is given by
\[ \tilde{a}(0)=\min_{\norm{z}=1}\sum_{k=1}^m |\ip{z}{f_k}|^4.  \]
\end{enumerate}
\label{t4.4}
\end{theorem}

The results presented so far show that both $\alpha$ and $\beta$ admit left inverses that are Lipschitz continuous on their domains of definition. 
One remaining problem is to know whether these left inverses can be extended to Lipschitz maps over the entire $\R^m$.
The following two results provide a positive answer (see \cite{BZ14,BZ15} for details).

\begin{theorem}[\cite{BZ15}] Assume $\fc\subset H=\C^n$ is a phase retrievable frame for $\C^n$. 
Let $\sqrt{A_0}$ be the lower Lipschitz constant of the map $\alpha:(\hat{H},D_2)\rightarrow(\R^m,\norm{\cdot}_2)$.
Then there is a Lipschitz map $\omega:(\R^m,\norm{\cdot}_2)\rightarrow(\hat{H},D_2)$ so that: (i) $\omega(\alpha(x))=x$ for all $x\in\hat{H}$, and (ii) its Lipschitz constant is $\Lip(\omega) \leq \frac{4+3\sqrt{2}}{\sqrt{A_0}}$.
\end{theorem}

\begin{theorem}[\cite{BZ14,BZ15a}] Assume $\fc\subset H=\C^n$ is a phase retrievable frame for $\C^n$. 
Let $\sqrt{a_0}$ be the lower Lipschitz constant of the map $\beta:(\hat{H},d_1)\rightarrow(\R^m,\norm{\cdot}_2)$.
Then there is a Lipschitz map $\psi:(\R^m,\norm{\cdot}_2)\rightarrow(\hat{H},d_1)$ so that: (i) $\psi(\beta(x))=x$ for all $x\in\hat{H}$, and (ii) its Lipschitz constant is $\Lip(\psi) \leq \frac{4+3\sqrt{2}}{\sqrt{a_0}}$.
\end{theorem}

{\bf Sketch of Proof}

Proofs of both results follow a similar strategy. First both metric spaces $(\hat{H},d_1)$ and $(\hat{H},D_2)$ are bi-Lipschitz isomorphic with $\SS^{1,0}$ via Lemma
\ref{l2.4}. Then one uses Kirszbraun's Theorem (see, e.g., \cite{BenLin00,HG13,WelWil75}) to obtain an isometric Lipschitz extension of the left inverse of $\alpha$ (or $\beta$) from its range to the entire $(\R^m,\norm{\cdot}_2)$ into $(\Sym(H),\norm{\cdot}_2)$. The final step is to construct a Lipschitz map $\pi:\Sym(H)\rightarrow \SS^{1,0}(H)$ so that $\pi(x^*)=xx^*$
 for every $x\in H$. This map is realized as $\pi(A)=(\lambda_1-\lambda_2)P_1$, where $\lambda_1\geq\lambda_2$ are the two largest eigenvalues of $A$, 
and $P_1$ is the principal eigenprojector. Using the integration contour from \cite{ZB06} and Weyl's inequalities (see III.2 in \cite{Bhatia}) the authors of \cite{BZ15}
obtained that $\pi$ is Lipschitz with $\Lip(\pi)\leq 3+2\sqrt{2}$ for $\pi:(\Sym(H),\norm{\cdot}_{2})\rightarrow(\SS^{1,0},\norm{\cdot}_{2})$. 

\subsection{Fisher Information Matrices and Cramer-Rao Lower Bounds}

Throughout this section assume $\fc=\{f_1,\ldots,f_m\}\subset H=\C^n$ is a phase retrievable frame for $V$, where $V\subseteq H$ is a real linear space, and $x\in V$. 

Consider two measurement processes. The first is the Additive White Gaussian Noise (AWGN) model \index{model! AWGN}
\begin{equation}
y_k  =|\ip{x}{f_k}|^2 + \nu_k~~,~~1\leq k\leq m, \label{eq:4.2.1} 
\end{equation}
where $(\nu_k)_{1\leq k\leq m}$ are independent and identically distributed (i.i.d.) realizations of a normal random variable of zero mean and variance $\sigma^2$.
The second process is a non-Additive White Gaussian Noise (nonAWGN) model \index{model! nonAWGN} where the noise is added prior to taking the absolute value:
\begin{equation}
\label{eq:NAWGNM}
y_k = |\ip{x}{f_k}+\mu_k|^2 ~~,~~1\leq k\leq m,
\end{equation}
where $(\mu_k)_{1\leq k\leq m}$ are i.i.d. realizations of a Gaussian complex process with zero mean and variance $\rho^2$.

First we present the {\em Fisher Information matrices} $\I$ \index{matrix! Fisher Information Matrix}\index{Fisher Information Matrix}for these two processes. The general definition
of the Fisher Information Matrix is (see \cite{Kay2010})
$$\I(x)=\E[(\nabla_x \log\,p(y;x))(\nabla_x \log\,p(y;x))^T].$$   
Following \cite{BCMN13a} and \cite{Bal13a} for the AWGN model (\ref{eq:4.2.1}) we obtain:
\begin{equation}
\label{eq:Fisher}
\I^{\AWGN}(x) = \frac{4}{\sigma^2}\Rsc(\xi) = \frac{4}{\sigma^2}\sum_{k=1}^m \Phi_k \xi\xi^T \Phi_k, 
\end{equation}
where $\xi=\j(x)\in\R^{2n}$. In general $\I(x)$ has rank at most $2n-1$ because $J\xi$ is always in its kernel.

 In \cite{Bal15} the Fisher
information matrix for the nonAWGN model (\ref{eq:NAWGNM}) is shown to have the following form:
\begin{eqnarray}
\I^{\nonAWGN}(x) & = &  \frac{4}{\rho^4} \sum_{k=1}^m \left( G_1\left( \frac{\ip{\Phi_k\xi}{\xi}}{\rho^2}\right)  -1 \right) \Phi_k\xi\xi^*\Phi_k \label{eq:I-G1} \nonumber\\
& = & \hspace{-3mm}\frac{4}{\rho^2} \sum_{k=1}^m G_2\left( \frac{\ip{\Phi_k\xi}{\xi}}{\rho^2}\right) \frac{1}{\ip{\Phi_k\xi}{\xi}}\Phi_k\xi\xi^*\Phi_k, \label{eq:I}
\end{eqnarray} 
where the two universal scalar functions $G_1,G_2:\R^+\rightarrow\R^+$ are given by
\begin{eqnarray}
\label{eq:G1}
 G_1(a) & = & \frac{e^{-a}}{a} \int_0^{\infty} \frac{I_1^2(2\sqrt{at})}{I_0(2\sqrt{at})} te^{-t} dt = \frac{e^{-a}}{8a^3} \int_0^{\infty} \frac{I_1^2(t)}{I_0(t)} t^3 e^{-\frac{t^2}{4a}} dt \\
G_2(a) & = & a(G_1(a)-1), \nonumber
\end{eqnarray}
where $I_0$ and $I_1$ are the modifed Bessel functions of the first kind and order 0 and 1, respectively.
Both Fisher information matrices have the same null space spanned by $J\xi$. 

Next we present a lower bound on the variance of any unbiased estimator for $x$. Let $z_0\in V$ be a fixed vector. Define
\begin{equation}
V_{z_0} = \{ x\in V~~,~~\ip{x}{z_0} > 0, \}
\end{equation}
where $\ip{\cdot}{\cdot}$ is the complex scalar product in $H$. Set $E_{z_0}=\span_{\R}(V_{z_0})$ the real vector space spanned by $V_{z_0}$. Note $E_{z_0}=\{x\in V~,~\imag(\ip{x}{z_0})=0\}$. 

To make (\ref{eq:4.2.1}) identifiable we select the representative $x\in V_{z_0}$ of the class $\hat{x}$. 
This is a mild condition since it only asks for the class $\hat{x}$ not to be orthogonal to $z_0$ with 
respect to the scalar product of $H$. 
An estimator $\omega:\R^m\rightarrow E_{z_0}$ is unbiased if $\E[\omega(\beta(x)+\nu)]=x$
for all $x\in V_{z_0}$. Here the expectation is taken with respect to the noise random variable.

A careful analysis (see \cite{Bal13a}) of the estimation process shows that the {\em Cramer-Rao Lower Bound (CRLB)} \index{bound! Cramer-Rao Lower Bound} 
for either measurement process (\ref{eq:4.2.1}) and (\ref{eq:NAWGNM})
is given by $(\PPRR \I(x)\PPRR)^{\dagger}$, where $\PPRR$ is the orthogonal projection onto $\Vc_{z_0}=\j(E_{z_0})$ in $\R^{2n}$, and upper script $\dagger$ denotes the Moore-Penrose pseudo-inverse. Here $\I(x)$ stands for the Fisher information matrix $\I^{\AWGN}(x)$ or $\I^{\nonAWGN}(x)$. 
Then the covariance of any unbiased estimator $\omega:\R^m\rightarrow E_{z_0}$ is bounded as follows:
\begin{equation}
\Cov[\omega] \geq \left(\PPRR\I(x)\PPRR\right)^{\dagger}.
\end{equation}
In particular, if one choses $z_0=x$ then $\PPRR$ becomes the orthogonal projection onto the range of $\I(x)$
 and $\left(\PPRR\I(x)\PPRR\right)^{\dagger}=\I(x)^{\dagger}$ (see (\ref{eq:CRLB22}) below).

 In the real case, $\fc\subset V=\R^n\subset\C^n$, the Fisher information matrices of the AWGN model (\ref{eq:4.2.1}) and of the nonAWGN model (\ref{eq:NAWGNM}) 
take the form
\[ \I^{\AWGN}(x) = \frac{4}{\sigma^2} \left[ \begin{array}{cc}
\mbox{$\Rbc(x)$} & 0 \\ 0 & 0 \end{array} \right]
~,~
\I^{\nonAWGN}(x) = \frac{4}{\sigma^2} \sum_{k=1}^m G_2(\frac{|\ip{x}{f_k}|^2}{\rho^2}) \left[ \begin{array}{cc}
\mbox{$f_kf_k^T$} & 0 \\ 0 & 0 \end{array} \right]. \]
Restricting to the real component of the estimator, the CRLB for the AWGN model (\ref{eq:4.2.1}) becomes
\[ \Cov[\omega^{\AWGN}] \geq \frac{\sigma^2}{4} \Rbc(x)^{-1} \]
whereas for the nonAWGN model (\ref{eq:NAWGNM}), the CRLB becomes
\[ \Cov[\omega^{\nonAWGN}] \geq \frac{\rho^2}{4} \left( \sum_{k=1}^m G_2(\frac{|\ip{x}{f_k}|^2}{\rho^2})f_kf_k^T \right)^{-1}. \]

In the complex case $\fc\subset V=H=\C^n$, $\PPRR=I_{2n}-J\psi_0\psi_0^TJ^T$ with $\psi_0=\j(z_0)$ and the CRLB for AWGN becomes
\[ \Cov[\omega] \geq \frac{\sigma^2}{4} ((I_{2n}-J\psi_0\psi_0^TJ^T)\Rsc(\xi)(I_{2n}-J\psi_0\psi_0^TJ^T))^{\dagger}. \]
Since $\fc$ is phase retrievable for $H$, by Theorem \ref{t4.4}(9) $\Rsc(\xi)$ satisfies the lower bound $\Rsc(\xi)\geq a_0 (\norm{\xi}^2 I_{2n}-J\xi\xi^TJ^T)$. A little bit of algebra shows


\begin{eqnarray*}
a_0 |\ip{x}{z_0}|^2 \Pi_{z_0} & = & a_0  |\ip{\xi}{\psi_0}|^2 \Pi_{z_0} \\
 & \leq &  a_0 |\ip{\xi}{\psi_0}|^2 (I_{2n}-J\psi_0\psi_0^TJ^T) + 
a_0 \left[ \norm{\xi-\ip{\xi}{\psi_0}\psi_0}^2 I_{2n} \right. \\
 & & \left. - J(\xi - \ip{\xi}{\psi_0}\psi_0)(\xi-\ip{\xi}{\psi_0}\psi_0)^TJ^T \right]  \\
 & = &   (I_{2n}-J\psi_0\psi_0^TJ^T)\Rsc(\xi)(I_{2n}-J\psi_0\psi_0^TJ^T). 
\end{eqnarray*}

In particular this inequality shows that, if an unbiased estimator $\omega^{0}:\R^m\rightarrow E_{z_0}$ for the AWGN model (\ref{eq:4.2.1})
achieves the CRLB then its covariance matrix is upper bounded by
\[ \Cov[\omega^{0}] \leq \frac{\sigma^2}{4a_0 |\ip{x}{z_0}|^2} \Pi_{z_0}. \]
This result was derived in \cite{Bal13a}.

Finally, if the global phase is provided by an oracle by correlating the estimated signal with the original signal $x$, then we can choose $z_0=x$ and $\Pi_{x}=I_{2n}-J\xi\xi^TJ^T$. 
But then
\[ \Pi_{x} \I^{\AWGN}(x)\Pi_x = \I^{\AWGN}(x)~~,~~\Pi_{x} \I^{\nonAWGN}(x)\Pi_x = \I^{\nonAWGN}(x) \]
which implies the CRLBs:
\begin{equation}
\label{eq:CRLB22}
\Cov[\omega^{\AWGN}] \geq \left(\I^{\AWGN}(x)\right)^{\dagger} ~~,~~\Cov[\omega^{\nonAWGN}] \geq \left(\I^{\nonAWGN}(x)\right)^{\dagger}.
\end{equation}

\section{Reconstruction Algorithms}

We present two types of reconstruction algorithms:
\begin{itemize}
\item Rank 1 tensor recovery: Linear Reconstruction, PhaseLift;
\item Iterative algorithms: Gerchberg-Saxton, Mean-Squares Optimization: \linebreak Wirtinger flow and IRLS.
\end{itemize}

The literature contains more algorithms than those presented here, see e.g. \cite{W12,ABFM12,Bal2010,FMNW13,fienup82}.

Throughout this section we assume $\fc$ is a phase retrievable frame for $H=\C^n$. We let $y=(y_k)_{1\leq k\leq m}$ denote the vector of measurements.
We analyze two cases: The noiseless case, when $y=\beta(x)$, and the additive noise case, when $y=\beta(x)+\nu$, where $\nu\in\R^m$ denotes the noise.

\subsection{Rank 1 Tensor Recovery}

The matrix recovery algorithms attempt to estimate the rank 1 matrix $X=xx^*$ from the measurements $y=(y_k)_{1\leq k\leq m}$.
We present two such algorithms: the linear reconstruction algorithm using lifting \index{lifting}, and PhaseLift. 
An extension of the linear reconstruction algorithm from a matrix to a higher order tensor setting is also included.

\subsubsection{Linear Reconstruction Using Lifting}
\mbox{}

{\bf (i) Order $2$ Tensor Embedding.}

Linear reconstruction works well when the frame has high redundancy. Specifically if $m\geq \dim_\R(\Sym(H))=n^2$ then, 
generically, the set of rank 1 operators $\{ F_k=f_kf_k^*~,~1\leq k\leq m\}$ is a frame for $\Sym(H)$. 
In this case the measurements are linear on the space of matrices:
\[ y_k=\ip{F_k}{X}_{HS} + \nu_k~~,~~1\leq k\leq m. \]
Let $\{ \tilde{F_k}~,~1\leq k\leq m\}$ denote the canonical dual frame to $\{ F_k~,~1\leq k\leq m\}$.
Then the minimum Frobenius norm estimate of $X$ is given by the linear formula ${X}_{\est} = \sum_{k=1}^m y_k \tilde{F_k}$.
The class $\hat{x}$ is recovered using the spectral decomposition of $X_{\est}$, 
$$X_{\est}=\sum_{j=1}^d \lambda_{r(j)}P_j,$$ 
where
$d$ denotes the number of distinct eigenvalues of $X_{\est}$, $\lambda_1\geq \lambda_2\geq\cdots\geq\lambda_n$ are the eigenvalues of $X_{\est}$, $P_j$ is the orthogonal projection onto the eigenspace associated to eigenvalue $\lambda_{r(j)}$, and $r(j+1)=r(j)+\rank(P_j)$,  with $r(1)=1$. 
The least squares estimator $x^{LS}$ of $\hat{x}$ from $X_{\est}$ minimizes $\norm{X_{\est}-xx^*}^2$. The solution is unique when the top eigenvalue of $X_{\est}$ is
simple, $\lambda_1>\lambda_2$. In this case let $P=ee^*$ for some unit norm vector $e$. The
least squares estimator $x^{LS}$ is given explicitly by
\begin{equation}
 x^{LS} = \left\{ \begin{array}{ccc}
\mbox{$\sqrt{\lambda_1}e$} & if & \mbox{$\lambda_1\geq 0$} \\
 0 & & otherwise \end{array} \right. .
\end{equation}
In the case $\lambda_1=\lambda_2>0$ there are infinitely many possible top eigenvectors. 
Unfortunately any such choice produces an estimator that is discontinuous as a function of $y$. 
On the other hand the following estimator
\begin{equation}
 x^{{\textnormal{Lip}}} = \left\{ \begin{array}{ccc}
\mbox{$\sqrt{\lambda_1-\lambda_2}e$} & if & \mbox{$\lambda_1 > \lambda_2$} \\
 0 & & otherwise \end{array} \right.
\end{equation}
 is an exact reconstruction scheme in the absence of noise and it is a Lipschitz continuous map with respect to the measurement vector $y$. Its Lipschitz constant
and its performance with respect to additive noise are described in \cite{BZ15}.
\vspace{10mm}

{\bf (ii) Higher Order Tensor Embeddings.}

The idea of lifting into the space $\SS^{1,0}$ of rank 1 matrices can be extended to spaces of higher order tensors (see \cite{Bal2009}). 
Fix an integer $d\geq 1$ and denote by $O_{n,d}=\{\gamma\in\N^d~,~1\leq \gamma(1)\leq\cdots\leq\gamma(d)\leq n\}$ the set of ordered 
$d$-tuples of positive integers up to $n$.
Let $\Pp_{d,d}(Z_1,\ldots,Z_n)$ denote the real linear space of symmetric homogeneous polynomials in $n$ variables $Z_1,\ldots,Z_n$ of degree $(d,d)$,
meaning that each monomial has degree $d$ in variables $Z_1,\ldots,Z_n$ and degree $d$ in conjugate variables $\overline{Z_1},\ldots,\overline{Z_n}$,
\[ P = \sum_{\gamma,\delta\in O_{n,d}} c_{\gamma,\delta} Z_{\gamma(1)}\cdots Z_{\gamma(d)}\overline{Z_{\delta(1)}\cdots Z_{\delta(d)}}
~~,~~ c_{\gamma,\delta}=\overline{c_{\delta,\gamma}}\in\C. \]
In the case $d=1$, $\Pp_{1,1}(Z_1,\ldots,Z_n)$ is a linear $\R$-space isomorphic to $\Sym(\C^n)$. In general, $\Pp_{d,d}(Z_1,\ldots,Z_n)$ is isomorphic to the $\R$-linear space
 of $(d,d)$-sesquilinear functionals over $\C^n$ (denoted by $\Lambda_{d,d}(\C^n)$ in \cite{Bal2009}). 
For a given ordered $d$-tuple $\gamma\in O_{n,d}$ we denote by $\Pi(\gamma)$ the collection of all permutations of $d$ elements that produce distinct $d$-tuples when applied to $\gamma$.
Let $d_1,d_2,\ldots,d_n$ denote, respectively, the number of repetitions of $1,2,\ldots,n$ in $\gamma$. Then the cardinal of $\Pi(\gamma)$ is given by the multinomial
formula $\Card(\Pi(\gamma))=\frac{d!}{d_1!\cdots d_n!}$. On $\Pp_{d,d}(Z_1,\ldots,Z_n)$ consider the sesquilinear scalar product $\ipp{\cdot}{\cdot}$ \index{product! $\ipp{\cdot}{\cdot}$} so that
$$\{ Z^{(\gamma,\delta)}:=\left( \Card(\Pi(\gamma))\Card(\Pi(\delta)) \right)^{1/2}Z_{\gamma(1)}\cdots Z_{\gamma(d)}\overline{Z_{\delta(1)}\cdots Z_{\delta(d)}} ~;~\gamma,\delta\in O_{n,d} \}$$ 
is an orthonormal basis. 
Let 
$$\kappa_{d,d}:\widehat{\C^n}\times\cdots\times\widehat{\C^n}\rightarrow \Pp_{d,d}(Z_1,\ldots,Z_n)~,~ \kappa_{d,d}(x^1,\ldots,x^d) = 
\prod_{k=1}^d|x^k_1Z_1+\cdots + x^k_nZ_n|^{2}.$$
Note $x\mapsto \kappa_{d,d}(x,x,\ldots,x)$ is an embedding of $\widehat{\C^n}$ into $\Pp_{d,d}(Z_1,\ldots,Z_n)$.
Let $P=\kappa_{d,d}(x,x,\ldots,x)$ and $Q_{k_1,\ldots,k_d}=\kappa_{d,d}(f_{k_1},\ldots,f_{k_d})$. Then a little algebra shows that
\begin{equation} \ipp{P}{Q_{k_1,\ldots,k_d}} = |\ip{x}{f_{k_1}}|^2\cdots |\ip{x}{f_{k_d}}|^2. \end{equation}
Now the phase retrieval problem can be restated as the problem of finding a homogeneous polynomial $P$ of rank 1 (that is, of the form $\kappa_{d,d}(x,\ldots,x)$) so that
\begin{equation}
\label{eq:Pabc}
\ipp{P}{Q_{k_1,\dots,k_d}} = y_{k_1}\cdots y_{k_d}~~,~~\forall (k_1,\ldots,k_d)\in O_{m,d}.
\end{equation}
The number of equations in (\ref{eq:Pabc}) is 
$$M_{m,d}=\Card(O_{m,d})=\left(\begin{array}{c} \mbox{$m+d-1$} \\ \mbox{$d$} \end{array}\right), $$ whereas the dimension of the
real linear space $\Pp_{d,d}(Z_1,\ldots,Z_n)$ is 
$$N_{n,d} = \dim_\R\PP_{d,d}(Z_1,\ldots,Z_n) = (\Card(O_{n,d}))^2 =\left(\begin{array}{c} \mbox{$n+d-1$} \\ \mbox{$d$} \end{array}\right)^2.$$
If $d$ is so that $M_{m,d}\geq N_{n,d}$ and the set of polynomials $\mathcal{Q}=\{Q_{\eps}~;~\eps\in O_{m,d}\}$ forms a frame for $\Pp_{d,d}(Z_1,\ldots,Z_n)$ then $P$ can be obtained by
 solving a linear system (albeit of dimension growing exponentially with $n$ and $m$). In particular if the set $\mathcal{Q}$ forms a $d$-design the reconstruction is particularly simple. The case $d=2$ has been also explored in \cite{BacEhl15}.
 In the absence of noise $P\in \kappa_{d,d}(\widehat{\C^n})$ and thus $\hat{x}$ is found by solving a factorization problem. In the presence of noise, $P$ is no longer of rank 1 and a different estimation procedure should be used. For instance one can find the ``closest" rank 1 homogeneous polynomial in $\Pp_{d,d}(Z_1,\ldots,Z_n)$ and invert $\kappa_{d,d}$ to estimate $x$.
\vspace{3mm}

\subsubsection{PhaseLift}

Consider the noiseless case $y=\beta(x)$. \index{PhaseLift}\index{algorithm! PhaseLift}
The main idea is embodied in the following feasibility problem:
\[ {\rm find}_{{\rm subject~ to:} \Ac(X)=y, X=X^*\geq 0, \rank(X)=1} X. \]
Except for the condition $\rank(X)=1$, the optimization problem would be convex. However the rank constraint destroys the convexity property.
Once a solution $X$ is found, the vector $x$ can be obtained by solving the factorization problem $X=xx^*$. 

The feasibility problem admits at most a unique solution and so does the following optimization problem:
\begin{equation}
\label{eq:opt0}
 \min_{\Ac(X)=y,X=X^*\geq 0} \rank(X), 
\end{equation}
which is still non-convex. The insight provided by matrix completion theory and exploited in \cite{CSV12,CESV12}
is to replace $\rank(X)$ by $\trace(X)$ which is convex. Thus one obtains
\begin{equation}
\label{eq:opt1}
{\rm (PhaseLift)}~~~~~~
\min_{\Ac(X)=y,X=X^*\geq 0} \trace(X),
\end{equation}
which is a convex optimization problem (a semi-definite program). In \cite{CL12} the authors proved that
for random frames, with high probability, the problem (\ref{eq:opt1}) has the same solution as the problem (\ref{eq:opt0}):
\begin{theorem} \label{t5.1}
Assume each vector $f_k$ is drawn independently from $\Nc(0,I_n/2)+i\Nc(0,I_n/2)$, or each vector is drawn independently from the uniform distribution on the complex sphere of radius $\sqrt{n}$. Then there are universal constants $c_0,c_1,\gamma>0$ so that for $m\geq c_0 n$, for every $x\in\C^n$ the problem (\ref{eq:opt1}) has
the same solution as (\ref{eq:opt0}) with probability at least $1-c_1 e^{-\gamma n}$. 
\end{theorem}
As explained in \cite{DemHan12} and \cite{CL12}, the minimization of trace is not necessary; in the absence of noise it reduces to a feasibility problem.

The PhaseLift algorithm is also robust to noise. Consider the measurement process
\[ y= \beta(x)+\nu, \]
for some $\nu\in\R^m$ noise vector.  Consider the following modified optimization problem:
\begin{equation}
\label{eq:opt2}
\min_{X=X^*\geq 0} \norm{\Ac(X)-y}_1.
\end{equation}
In \cite{CL12} the following result has been shown:
\begin{theorem}\label{t5.2}
Consider the same stochastic process for a random frame $\fc$. There is a universal constant $C_0>0$ so that
for all $x\in C^n$ the solution to (\ref{eq:opt2}) obeys
\[ \norm{X - xx^*}_2 \leq C_0 \frac{\norm{\nu}_1}{m}. \]
For the Gaussian model this holds with the same probability as in the noiseless case, whereas the probability of failure is exponentially small in $n$ in the uniform model. The principal eigenvector $x^0$ of $X$ (normalized by the square root of the principal eigenvalue) obeys
\[ D_2(x^0,x) \leq C_0 \min(\norm{x}_2,\frac{\norm{\nu}_1}{m\norm{x}_2}) .\] 
\end{theorem}

\subsection{Iterative Algorithms}

We present two classes of iterative algorithms: the Gerchberg-Saxton algorithm and mean-squares minimization algorithms.

\subsubsection{The Gerchberg-Saxton Algorithm}

Let $c=(c_k)_{1\leq k\leq m}\in\C^m$ denote a sequence of the frame coefficients $c_k=\ip{x}{f_k}$. \index{algorithm! Gerchberg-Saxton}
Let $E = \{ (\ip{x}{f_k})_{1\leq k\leq m}~,~x\in \C^n\}$ denote the range of frame coefficients. Assume the measurements are all nonnegative,
 $y_k\geq 0$ for all $k$ (otherwise rectify at $0$). Denote by $\{\tilde{f_k}~,~1\leq k\leq m\}$ the canonical dual frame of $\{f_1,\ldots,f_m\}$. 
The Gerchberg-Saxton algorithm first introduced in \cite{GS72} iterates between two sets of constraints: $|c_k|=\sqrt{y_k}$, and $c\in E$.
Let $x^{0}\in H$ be an initialization and set $t=0$.
The algorithm repeats the following steps:
\begin{enumerate}
\item Linear Analysis: $c_k=\ip{x^{t}}{f_k}$, $1\leq k\leq m$;
\item Magnitude Adjustment: $d_k= \sqrt{y_k}c_k/|c_k|$, $1\leq k\leq m$;
\item Linear Synthesis: $x^{t+1}=\sum_{k=1}^m d_k \tilde{f_k}$;
\item Increment $t=t+1$;
\end{enumerate}
until a stopping criterion is achieved. The main advantage of this algorithm is its simplicity. It can easily incorporate additional constraints on $x$ or its transform $c$.
Unfortunately it suffers of a couple of disadvantages. Namely, the convergence is not guaranteed, and furthermore, when it converges it only converges to a local minimum. 
However, despite these shortcomings, the algorithm performs relatively well when $x$ is highly constrained, for instance when all of its entries are non-negative 
(see e.g. \cite{fienup82}).

\subsubsection{Mean-Squares Minimization Algorithms}

Consider again the measurement process
\[ y_k = |\ip{x}{f_k}|^2 + \nu_k~~,~~1\leq k\leq m. \]
The Least-Squares criterion
\[ \min_{x\in\C^n} \sum_{k=1}^m ||\ip{x}{f_k}|^2-y_k|^2 \]
can be understood as the Maximum Likelihood Estimator (MLE) when the additive noise vector $\nu\in\R^m$ is
normal distributed with zero mean and covariance $\sigma^2 I_m$. However the optimization problem 
is not convex and has many local minima. 

We present two algorithms that minimize the mean-squares error: the Wirtinger flow and the Iterative Regularized Least-Squares.

{\bf (i) Gradient Descent Using the Wirtinger Flow}

This algorithm has been introduced in \cite{CLS14}. The idea is to follow the gradient descent for the criterion \index{algorithm! Wirtinger flow}
\begin{equation}
f(x) = \sum_{k=1}^m |y_k - |\ip{x}{f_k}|^2|^2.
\end{equation}
The initialization is performed using the spectral method. Specifically:

{\bf Step 1. Initialization.} Compute the principal eigenvector of \linebreak $R_y=\sum_{k=1}^m y_k f_kf_k^*$ using, e.g., the power method. Let $(e_1,a_1)$ be the eigen-pair with $e_1\in\C^n$, $\norm{e_1}=1$, and $a_{1}\in\R$. 
Initialize:
\begin{equation}
\label{eq:Wstep1}
x^{0}  = \sqrt{n\frac{\sum_{k=1}^m y_k}{\sum_{k=1}^m \norm{f_k}^2}}\, e_1 ~,~t  = 0.
\end{equation}

{\bf Step 2. Iteration.} Repeat:

2.1 Gradient descent:
\begin{equation}
x^{t+1} = x^{t} - \frac{\mu_{t+1}}{\norm{x^{0}}^2} \left( \frac{1}{m} \sum_{k=1}^m (|\ip{x^{t}}{f_k}|^2-y_k)\ip{x^{t}}{f_k}f_k \right).
\end{equation}

2.2 Update
\begin{equation}
\mu^{t+1} = \min(\mu_{\textnormal{max}},1-e^{-\tau/\tau_0}).
\end{equation}

{\bf Step 3. Stopping.} Stop after a fixed number of iterations or an error criterion is achieved.

The authors of \cite{CLS14} showed that this algorithm converges with high probability to the exact solution in the absence of noise:
\begin{theorem}[\cite{CLS14}] Let $x\in\C^n$ and $y=\beta(x)$ with $m\geq c_0 n\log n$, where $c_0$ is a sufficiently large constant. Then
the Wirtinger flow initial estimate $x^0$, normalized to have squared Euclidean norm equal to $\frac{1}{m}\sum_{k=1}^m y_k$, obeys $D_2(x^0,x)\leq \frac{1}{8}\norm{x}$
 with probability at least $1-10e^{-\gamma n}-\frac{8}{n^2}$, where $\gamma$ is a fixed positive numerical constant. 
Further, take a constant learning parameter sequence,
 $\mu^{t}=\mu$ for all $t\geq 1$ and assume $\mu\leq \frac{c_1}{n}$ for some fixed numerical constant $c_1$. Then there is an event of probability at least $1-13 e^{-\gamma n}-m e^{-1.5 m}-\frac{8}{n^2}$, such that on this event, starting from any initial solution $x^0$ obeying $D_2(x^0,x)\leq \frac{1}{8}\norm{x}$, we have
\begin{equation}
D_2(x^t,x) \leq \frac{1}{8}\left(1-\frac{\mu}{4}\right)^{t/2} \norm{x}.
\end{equation}
\end{theorem}

{\bf (ii) The Iterative Regularized Least-Squares Algorithm}

The iterative algorithm described next tries to find the global minimum using a regularization term.
Consider the following optimization criterion:
\begin{eqnarray}
\label{eq:J}
 J(u,v;\lambda,\mu) & = & \sum_{k=1}^m \left|\frac{1}{2}(\ip{u}{f_k}\ip{f_k}{v}+\ip{v}{f_k}\ip{f_k}{u}) - y_k\right|^2 \\
 & & 
+\lambda \norm{u}_2^2 + \mu \norm{u-v}_2^2 + \lambda \norm{v}_2^2. \nonumber
\end{eqnarray}
The {\em Iterative Regularized Least-Squares (IRLS)} \index{algorithm! IRLS} algorithm presented in \cite{Bal13a} works as follows.

Fix a stopping criterion, such as a tolerance $\eps$, a desired level of signal-to-noise-ratio $\snr$, or a 
maximum number of steps $T$. 
Fix an initialization parameter $\rho\in (0,1)$, a learning rate $\gamma\in(0,1)$, and a saturation parameter $\mu_{\textnormal{min}}>0$. 

{\bf Step 1. Initialization.} Compute the principal eigenvector of \linebreak $R_y=\sum_{k=1}^m y_k f_kf_k^*$ using e.g. the power method. Let $(e_1,a_1)$ be the eigen-pair with $e_1\in\C^n$, $\norm{e_1}=1$, and $a_{1}\in\R$.
If $a_{1}\leq 0$ then set $x=0$ and exit. Otherwise initialize:
\begin{equation}
\label{eq:step1}
x^{0}  = \sqrt{\frac{(1-\rho)a_1}{\sum_{k=1}^m |\ip{e_1}{f_k}|^4}}\, e_1~,~\lambda_0  =  \rho a_1~,~\mu_0  =  \rho a_1 ~,~t  = 0.
\end{equation}

{\bf Step 2. Iteration.} Perform:

2.1 Solve the least-squares problem:
\[ x^{t+1} = {\rm argmin}_u J(u,x^t;\lambda_t,\mu_t) \]
using the conjugate gradient method.

2.2 Update:
\[ \lambda_{t+1} = \gamma \lambda_t~,~\mu_t = \max(\gamma\mu_t,\mu_{\textnormal{min}})~~,~~t=t+1. \]

{\bf Step 3. Stopping.} Repeat Step 2 until
\begin{itemize}
\item The error criterion is achieved: $J(x^t,x^t;0,0) <\eps$; 
\item The desired signal-to-noise-ratio is reached: $\frac{\norm{x^t}^2}{J(x^t,x^t;0,0)} > \snr$; or
\item The maximum number of iterations is reached: $t>T$.
\end{itemize}

The final estimate can be $x^T$ or the best estimate obtained in the iteration path: $x^{\est} = x^{t_0}$ where $t_0 ={\rm argmin}_t J(x^t,x^t;0,0)$. 

The initialization (\ref{eq:step1}) is performed for the following reason. Consider the modified criterion:
\begin{eqnarray*}
 H(x;\lambda)  & = & J(x,x;\lambda,0) = \norm{\beta(x)-y}_2^2 + 2\lambda \norm{x}_2^2  \\
 & = & \sum_{k=1}^m |\ip{x}{f_k}|^4 + 2\ip{(\lambda I_n - R_y)x}{x} + \norm{y}_2^2. 
\end{eqnarray*}
In general this function is not convex in $x$, except for large values of $\lambda$. Specifically for $\lambda>a_1$,
 the largest eigenvalue of $R_y$, $x\mapsto H(x;\lambda)$ is convex and has a unique global minimum at $x=0$. 
For $a_1-\eps<\lambda<a_1$ the criterion is no longer convex, but the global minimum stays in a neighborhood of the origin. Neglecting the $4^{th}$ order terms, the critical points are given by the eigenvectors of $R_y$. Choosing 
$\lambda=\rho a_1$ and $x=s e_1$, the optimal value of $s$ for $s\mapsto H(se_1;\rho a_1)$ is given in (\ref{eq:step1}). 

The path of iterates $(x^t)_{t\geq 0}$ can be thought of as trying to approximate the measured vector $y$ 
with $\Ac(\outp{x^{t-1}}{x^t})$, where $\Ac$ is defined in (\ref{eq:Ac}). 
The parameter $\mu$ penalizes the unique negative eigenvalue of $\outp{x^{t-1}}{x^t}$; the larger the value of $\mu_t$ the smaller the iteration step $\norm{x^{t+1}-x^t}$ and the smaller the deviation of the matrix $\outp{x^{t+1}}{x^t}$ from a rank 1 matrix; the smaller the parameter $\mu_t$ the larger in magnitude the negative eigenvalue of $\outp{x^{t+1}}{x^t}$. This fact explains why in the noisy case the iterates first decrease the matching error $\norm{\Ac(x^t (x^t)^*-y}_2$ to some $t_0$, and then they may start to 
increase this error; instead the rank 2 self-adjoint operator $T=\outp{x^{t+1}}{x^t}$ always decreases 
the matching error $\norm{\Ac(T)-y}_2$.

At  any point on the path, if the value of criterion $J$ is smaller than the value reached at the target vector $x$, then the algorithm is guaranteed to converge near $x$. 
Specifically in \cite{Bal13a} the following result has been proved:
\begin{theorem}[\cite{Bal13a} Theorem 5.6] Fix $0\neq z_0\in\C^n$.
Assume the frame $\fc$ is so that $\ker\Ac\cap\SS^{2,1} =\{0\}$. Then there is a constant $A_3>0$ that depends on $\fc$ so that for every $x\in \C^n$ with $\ip{x}{z_0}>0$ and $\nu\in\C^n$ that produce $y=\beta(x)+\nu$ 
if there are $u,v\in\C^n$ so that $J(u,v;\lambda,\mu)<J(x,x;\lambda,\mu)$ then
\begin{equation}
\norm{\outp{u}{v} - xx^*}_1 \leq \frac{4\lambda}{A_3} + \frac{2\norm{\nu}_2}{\sqrt{A_3}}.
\end{equation}
Moreover, let $\outp{u}{v}=a_{1}e_1e_1^* + a_{2}e_2e_2^*$ be its spectral factorization with $a_{1}\geq 0\geq a_{2}$ and $\norm{e_1}=\norm{e_2}=1$. Set $\tilde{x}=\sqrt{a_{1}}e_1$. Then
\begin{equation}
 D_2(x,\tilde{x})^2 \leq \frac{4\lambda}{A_3} + \frac{2\norm{\nu}_2}{\sqrt{A_3}} + \frac{\norm{\nu}^2_2}{4\mu} + \frac{\lambda\norm{x}_2^2}{2\mu}.
\end{equation}
\end{theorem}
The kernel requirement on $\Ac$ is satisfied for generic frames when $m\geq 6n$. In particular this condition 
requires the frame $\fc$ is phase retrievable for $\C^n$.

\bibliographystyle{amsplain}

\end{document}